\documentclass[12pt,a4paper]{article}
\usepackage{amsfonts}
\usepackage{amsmath}
\usepackage{mathrsfs}
\usepackage{graphicx}
\usepackage{amsmath,amsthm,amssymb,amscd}
\newtheorem{Thm}{Theorem}[section]
\newtheorem{Lem}[Thm]{Lemma}

\newtheorem{Def}[Thm] {Definition}
\newtheorem{Ex}  [Thm]{Example}
\newtheorem{Cor}[Thm]{Corollary}
\theoremstyle{remark}
\newtheorem{Rem} [Thm]{Remark}

\theoremstyle{claim}

\newtheorem{Que}[Thm]{Question}

\DeclareMathOperator{\Diff}{Diff}

\begin{document}

\begin{center}
{\Large \bf  Minimal Diffeomorphisms cannot satisfy Generalized Dominated Splitting
}\\
\smallskip
\end{center}

\bigskip

\begin{center}
Xueting Tian $^\dagger$
\end{center}
\begin{center}
Academy of Mathematics and Systems Science, Chinese Academy of Sciences, Beijing 100190, China\\
\end{center}
\begin{center}
E-mail: tianxt@amss.ac.cn; txt@pku.edu.cn
\end{center}
\bigskip

 \footnotetext{$^\dagger$ Tian is  supported by
CAPES.}
  \footnotetext {AMS Review:  37D30, 37C40, 37B20.}

\smallskip
\begin{abstract}

We introduce a new notion called generalized dominated splitting which is  weaker than classical dominated splitting. We use this notion to generalize a result of Zhang\cite{Zh}: every diffeomorphism with nontrivial global generalized dominated splitting can not be minimal.
\end{abstract}

\section{Introduction}\label{Introduction}
Minimality is an important concept in the study of dynamical systems. It is interesting to study some
nature structure of the system that incompatible with minimality. In 1980's Herman\cite{Herman}  constructed a (family of) $C^1$ diffeomorphism on a compact 4-dimensional manifold
that is minimal but has positive topological entropy simultaneously and Rees\cite{Rees} constructed a minimal
homeomorphism with positive topological entropy on 2-torus. So positive entropy
is insufficient to guarantee the non-minimality. In \cite{Mane} Ma\~{n}\'{e} gave an
argument to locate some nonrecurrent point if the map admits some invariant expanding
foliation (also see \cite{BHH}). In particular this argument implies that a partially hyperbolic
diffeomorphism always has some nonrecurrent point and hence can not be minimal. Recently in \cite{Zh} Zhang showed that a global dominated splitting is sufficient to exclude the minimality
of the system. In present paper we mainly want to generalize the result  \cite{Zh} to a more general assumption called generalized dominated splitting.


Let $M$ be a compact $D$-dimensional smooth Riemannian manifold and
let $d$ denote the distance induced by the Riemannian metric.
Denote the tangent bundle of $M$  by $TM$ and denote by $\Diff^1(M)$ the space of $C^{1}$
diffeomorphisms of $M$. Denote  the maximal norm of  a linear map $A$ by $\|A\|$ and denote the minimal norm of an invertible linear map $A$ by
$m(A):=\|A^{-1}\|^{-1}$. Now we introduce our new notion called generalized dominated splitting.

\begin{Def}\label{def:gen-dom}Let $f\in \Diff^1(M)$ and let $S\in \mathbb{N}, \lambda>0$.
Given an $f-$invariant compact set $\Delta$, we say a $Df-$invariant splitting $T_{\Delta}M=E\oplus F$ on
$\Delta$ to be a {\bf generalized dominated splitting} on $\Delta$ (or simply GDS), if \\
(1). $T_{\Delta}M=E\oplus F$  is continuous  on
$\Delta$;\\
(2).
$\frac
{\|Df^{kS}|_{E(x)}\|}{m(Df^{kS}|_{F(x)})}\leq \lambda,\,\,\forall x \in
\Delta,\,\,\forall k\in \mathbb{N};$ \\
(3). there exists $x_0\in \Delta$,
$\frac
{\|Df^S|_{E(x_0)}\|}{m(Df^S|_{F(x_0)})}< \lambda^{-1}.$
\end{Def}

Note that $\|AB\|\leq \|A\| \|B\|$, $m(AB)\geq m(A)m(B)$, $\|AB\|\geq \|A\| m(B)$ and $m(AB)\leq m(A)\|B\|$. Then
\begin{eqnarray}\label{Eq:NormOperater}
\frac{\|Df^{[\frac nS]S}|_{E(x_0)}\|}{m(Df^{[\frac nS]S}|_{F(x_0)})}\times C^{2S} \geq \frac{\|Df^n|_{E(x_0)}\|}{m(Df^n|_{F(x_0)})}\geq \frac{\|Df^{[\frac nS]S}|_{E(x_0)}\|}{m(Df^{[\frac nS]S}|_{F(x_0)})}\times C^{-2S}
\end{eqnarray}
where $C=sup_{x\in M}\max\{\|Df(x)\|,\|Df^{-1}(x)\|\}$.

\begin{Rem}\label{Rem:DefGDS}
By (\ref{Eq:NormOperater}) the second condition in Definition \ref{def:gen-dom} implies that for all $x\in\Delta$, $$\limsup_{n\rightarrow +\infty}\frac1n\log\frac{\|Df^{n}|_{E(x)}\|}{m(Df^{n}|_{F(x)})}=\frac1S\limsup_{n\rightarrow +\infty}\frac1n\log\frac{\|Df^{nS}|_{E(x)}\|}{m(Df^{nS}|_{F(x)})}\leq \frac1S\min\{\log \lambda, 0\}\leq0.$$
\end{Rem}

Recall the definition of  classical dominated splitting. Let $\Delta$ be an $f-$invariant compact set. A $Df-$invariant splitting $T_{\Delta}M=E\oplus F$ on
$\Delta$ is called $(S, \lambda)
$-dominated on $\Delta$ (or simply dominated), if there exist two
constants $S\in \mathbb{Z}^+$ and $0<\lambda<1$ such that
$$\frac
{\|Df^S|_{E(x)}\|}{m(Df^S|_{F(x)})}\leq\lambda,\,\,\forall x \in
\Delta.$$  Note that dominated splitting is always continuous  (see \cite{BLV}), $\lambda<1<\lambda^{-1}$ and $$\frac{\|Df^{kS}|_{E(x)}\|}{m(Df^{kS}|_{F(x)})}\leq\prod_{i=0}^{k-1}\frac{\|Df^S|_{E(f^i(x))}\|}{m(Df^S|_{F(f^i(x))})}\leq\lambda^k\leq \lambda.$$ So   any dominated splitting satisfies Definition \ref{def:gen-dom}(Moreover, we will give an example below which does not have global dominated splitting but admits a global GDS).

\medskip

Let $\Delta$ be an $f-$invariant compact set. If $\varnothing\neq\Delta\varsubsetneqq M$, it is clear that $f$ can not be minimal since the closure of every orbit in $\Delta$ is still contained in $\Delta$.
A $Df-$invariant splitting $T_{\Delta}M=E\oplus F$ on
$\Delta$  is nontrivial if $dim(E)\cdot dim(F)\neq 0.$ And we say a $Df-$invariant  splitting $T_{\Delta}M=E\oplus F$ to be global, if $\Delta=M.$ Now we  state our main theorem for considering  systems with global GDS.

\begin{Thm}\label{MainThm}
Let $f\in \Diff^1(M)$. If there is a nontrivial global GDS $TM=E\oplus F$, then $f$ can not be minimal.
\end{Thm}



\begin{Rem}\label{Rem:GDS}
There exists  some minimal system $f\in \Diff^1(M)$ such that its
nontrivial global invariant splitting $TM=E\oplus F$ only satisfies the
first two conditions in Definition \ref{def:gen-dom}. More precisely, for
$a\in \mathbb{R}$, define the corresponding rotation $R_a: \mathbb{S}^1\rightarrow
\mathbb{S}^1,\,x \mapsto x+a(mod 1)$. Clearly for the product  system $f_{a_1,a_2,\cdots,a_n}:=R_{a_1}\times R_{a_2}\times \cdots \times R_{a_n}:\mathbb{T}^n\rightarrow \mathbb{T}^n$($a_1,a_2,\cdots,a_n\in\mathbb{R}$),
its nontrivial global invariant splitting $TM=E\oplus F$  satisfies the first two conditions
in Definition \ref{def:gen-dom} for any $S\in\mathbb{N},\lambda\geq1$ but fails the
third condition for all $S\in\mathbb{N},\lambda\geq1$ ($\frac
{\|Df^{k}|_{E(x)}\|}{m(Df^{k}|_{F(x)})}\equiv 1,\,\,\forall x \in
\Delta,\,\,\forall k\in \mathbb{N}$).
It is well-known that  $1,a_1,a_2,\cdots,a_n\in\mathbb{R}$ are  rationally independent if and only if 
 the product  system $f_{a_1,a_2,\cdots,a_n}=R_{a_1}\times R_{a_2}\times \cdots \times R_{a_n}$ 
 is minimal(and ergodic with respect to Lebesgue measure).
  This shows that both minimal and non-minimal $C^{\infty}$ diffeomorphisms admit to have nontrivial global
  invariant splitting $TM=E\oplus F$  satisfying the first two conditions in Definition \ref{def:gen-dom}.
\end{Rem}


\medskip

If  a global $Df-$invariant  splitting $T_{\Delta}M=E\oplus F$ is not GDS but the first two conditions in Definition \ref{def:gen-dom} still hold for some $S,\lambda$, then $\lambda\geq1$ and  $\lambda^{-1}\leq\frac
{\|Df^{kS}|_{E(x)}\|}{m(Df^{kS}|_{F(x)})}\leq \lambda,\,\,\forall x \in
\Delta,\,\,\forall k\in \mathbb{N}.$ Otherwise,  there exists $x_0\in \Delta$ and $k_0$,
$\frac
{\|Df^{k_0S}|_{E(x_0)}\|}{m(Df^{k_0S}|_{F(x_0)})}< \lambda^{-1},$ then $T_{\Delta}M=E\oplus F$ is GDS for $N:=k_0S$ and $\lambda.$ Furthermore, from (\ref{Eq:NormOperater}) $$(C^{2S}\lambda)^{-1}\leq\frac
{\|Df^{n}|_{E(x)}\|}{m(Df^{n}|_{F(x)})}\leq C^{2S} \lambda,\,\,\forall x \in
\Delta,\,\,\forall n\in \mathbb{N}.$$  In particular, $\lim_{n\rightarrow +\infty}\frac1n\log\frac{\|Df^{n}|_{E(x)}\|}{m(Df^{n}|_{F(x)})}=0,\,\,\forall x \in
\Delta.$
From these analysis and Remark \ref{Rem:DefGDS}, the third condition in Definition \ref{def:gen-dom} can be deduced once for some $x$, $$\liminf_{n\rightarrow +\infty}\frac1n\log\frac{\|Df^{n}|_{E(x)}\|}{m(Df^{n}|_{F(x)})}<0.$$
Hence, we state such a corollary of Theorem \ref{MainThm} as follows.
\begin{Cor}\label{Cor:MainThm}
Let $f\in \Diff^1(M)$. If there is a nontrivial global  $Df-$invariant splitting $TM=E\oplus F$ satisfying the first two conditions in Definition \ref{def:gen-dom} and there exists $x_0$, $$\liminf_{n\rightarrow +\infty}\frac1n\log\frac{\|Df^n|_{E(x_0)}\|}{m(Df^n|_{F(x_0)})}<0,$$ then $TM=E\oplus F$ is a nontrivial global GDS and thus $f$ can not be minimal.
\end{Cor}

This corollary can be as a sufficient condition to obtain GDS. By Remark \ref{Rem:DefGDS}, for a global $(S,\lambda)$-dominated splitting, every point $x$ satisfies $$\limsup_{n\rightarrow +\infty}\frac1n\log\frac{\|Df^n|_{E(x)}\|}{ m(Df^n|_{F(x)})}\leq{\frac 1S} \log\lambda<0.$$
This implies for any global dominated splitting, every point  satisfies the assumption of Corollary \ref{Cor:MainThm} and the supreme limit can be uniformly less than  $0$. But Corollary \ref{Cor:MainThm} assumes only one such point and uses inferior limit. Thus the assumption in Corollary \ref{Cor:MainThm} is still weaker than dominated splitting(for instance, see Example \ref{Ex:GDS-DS} below).

\begin{Rem}\label{Rem:Cor-GDS}For  any {\it surface} diffeomorphism $f$  with positive topological entropy,
if there is a nontrivial global  $Df-$invariant splitting $TM=E\oplus F$ satisfying the first
two conditions in Definition \ref{def:gen-dom}, then $f$ satisfies Corollary \ref{Cor:MainThm}
and thus $f$ can not be minimal. More precisely, by  Variational Principle(\cite{Walter}), for
any diffeomorphism with positive entropy, there exists an ergodic measure $\mu$ with positive entropy
and thus by Rulle's inequality(\cite{Ru}) $\mu$ has both negative Lyapunov exponent($\chi_1<0$) and positive
 Lyapunov exponent($\chi_2>0$) simultaneously. Note that $dim(E)= dim(F)=1,$ then  for $\mu$ a.e.
  $x$, $$\lim_{n\rightarrow +\infty}\frac1n\log\frac{\|Df^{n}|_{E(x)}\|}{m(Df^{n}|_{F(x)})}=\chi_1-\chi_2\leq -2h_{\mu}(f)<0.$$
  In particular, we recall that  if $f$ is a  $C^{1+\alpha}$ {\it surface} diffeomorphism
  with positive topological entropy, then $f$ always has periodic point by classical Pesin theory\cite{Katok2} and thus can not be minimal.

\end{Rem}


\section{Proof of Theorem \ref{MainThm}} \label{Main Proof}

Before proving  Theorem \ref{MainThm} we need a result of \cite{Herz} .

\begin{Lem}\label{Lem:1}(Proposition 3.4 in \cite{Herz}) Let $f : X \rightarrow X$ be a continuous map of a compact metric space.
Let $a_n : X \rightarrow R, n \geq 0,$ be a sequence of continuous functions such that
\begin{eqnarray}\label{Lem1Eq1}
a_{n+k}(x) \leq a_n(f^k(x)) + a_k(x) \text{ for every  } x \in X, n, k \geq 0.
\end{eqnarray}
and such that there is a sequence of continuous functions $b_n : X \rightarrow R, n \geq 0,$ satisfying
\begin{eqnarray}\label{Lem1Eq2}
 a_n(x) \leq a_n(f^k(x)) + a_k(x) + b_k(f^n(x)) \text{ for every }  x \in X, n, k \geq 0.
\end{eqnarray}
If $$\inf\frac1n\int_Xa_n(x) d\mu<0$$ for every ergodic f-invariant measure, then there is $N > 0$
such that $a_N(x) < 0$ for every $x \in X.$

\end{Lem}

\bigskip

{\bf Proof of Theorem \ref{MainThm}}
If $\lambda<1$, then the nontrivial global GDS is a nontrivial global dominated splitting. By the result of \cite{Zh}, $f$ can not be minimal. Now we assume that $\lambda\geq 1$ and we will give a proof by contradiction. Suppose $f$ is minimal, then the nontrivial global GDS can not be a nontrivial global dominated splitting from \cite{Zh}. To get a contradiction for this case, we only need to prove that the nontrivial global GDS is nontrivial global dominated splitting.

 Define for $\epsilon>0$ $$A_\epsilon:=\{z\in M\,\,\,|\,\,\,\frac
{\|Df^S|_{E(z)}\|}{m(Df^S|_{F(z)})}< -\epsilon+\lambda^{-1} \}$$ and set $A=\bigcup_{\epsilon>0} A_l.$  Note that every set $A_\epsilon$ is open set and $A$ is also open since the splitting $TM=E\oplus F$ is continuous. Clearly by assumption the point $x_0$ must be in $ A$.  Take $0<\epsilon<\lambda^{-1}$ small enough such that $x_0\in A_\epsilon$ so that $A_\epsilon\neq\varnothing$.

Since we assume that $f$ is minimal, then for every invariant measure $\mu,$ its support $supp(\mu)$ must coincide with the whole manifold $M$. Otherwise, if for some $\mu$, $supp(\mu)\subsetneqq M.$ Then every point $x\in supp(\mu)$, the closure of its orbit is contained in $supp(\mu)\subsetneqq M$ since $supp(\mu)$ is always compact and invariant. This contradicts that $f$ is minimal.
So for any nonempty open set, it always has positive measure for any invariant measure. In particular,  $\mu(A_\epsilon)>0$ holds for any invariant measure $\mu.$

Define functions for  $x\in M$ $$a_n(x):=\log\frac{\|Df^n|_{E(x)}\|}{ m(Df^n|_{F(x)})}.$$  Recall that $\|AB\|\leq \|A\| \|B\|$ and $m(AB)\geq m(A)m(B)$. Then it is easy to see that $a_n$ satisfy (\ref{Lem1Eq1}) of Lemma \ref{Lem:1}. Taking into account (\ref{Lem1Eq1}) we see that (\ref{Lem1Eq2}) holds once $a_n(x) \leq a_{n+k}(x) + b_k(f^n(x)).$ This is easily verified
for  $b_k(x):=\log\frac{\|(Df^{k}|_{E(x)})^{-1}\|}{ m((Df^{k}|_{F(x)})^{-1})}$
since $$\frac{\|Df^{n}|_{E(x)}\|}{ m(Df^{n}|_{F(x)})}\leq \frac{\|Df^{n+k}|_{E(x)}\|}{ m(Df^{n+k}|_{F(x)})}\times\frac{\|(Df^{k}|_{E(f^n(x))})^{-1}\|}{ m((Df^{k}|_{F(f^n(x))})^{-1})}.$$
Recall that  $TM=E\oplus F$ is a continuous splitting. So $a_n(x),b_n(x)$ are continuous functions.
Then all assumptions of Lemma \ref{Lem:1} are satisfied once  $\inf\frac1n\int_M a_n(x) d\mu<0$ holds for ergodic invariant measure $\mu$.

Let $\mu$ be an $f$ ergodic invariant measure.  By Subadditive Ergodic Theorem(see \cite{Walter}), and the ergodicity of $\mu$, the limit function $$a(x):=\lim_{n\rightarrow +\infty}\frac1n a_n(x)$$ is well-defined, $f$-invariant and can be a constant function for $\mu$ a.e $x$. Now we prove for  $\mu$ a.e $x$, $a(x)<0$ which implies $\inf\frac1n\int_M a_n(x) d\mu<0$.  Let $\Phi:=\cup_{n\in \mathbb{Z}}f^n(A_\epsilon)$. Clearly it is $f$-invariant and from ergodicity of $\mu$ we have that $\mu(\Phi)=1$(In fact, $\Phi=M$ since $M\setminus\Phi$ is $f$-invariant and closed but $f$ is minimal). So we only need  to prove $a(x)<0$ for $\mu$ a.e $x\in A_\epsilon$ since $a(x)$ is $f$ invariant.
Define $c_n(x):=a_{nS}(x)$, then $$c(x):=\lim_{n\rightarrow +\infty} \frac1n c_n(x)=S\lim_{n\rightarrow +\infty}\frac1{nS} a_{nS}(x)=Sa(x).$$ So we only need to prove $c(x)<0$ for $\mu$ a.e. $x\in A_\epsilon$.

By Birkhoff Ergodic Theorem,  the limit function $$\chi^*_{A_\epsilon}(x):=\lim_{n\rightarrow+\infty}\frac1n\sum_{j=0}^{n-1}\chi_{A_\epsilon}(f^{jS}(x))$$ exists for $\mu$ a.e. $x$ and $f^S-$invariant. Moreover, $\int\chi^*_{A_\epsilon}(x)d\mu=\int\chi_{A_\epsilon}(x)d\mu=\mu({A_\epsilon}).$ If $\mu$ is $f^S$ ergodic, it is obvious  since $\chi^*_{A_\epsilon}(x)\equiv\int\chi_{A_\epsilon}(x)d\mu=\mu({A_\epsilon})>0$ holds for $\mu$ a.e. $x\in M$. But we do not know whether $\mu$ is $f^S$ ergodic or not. We only know that it is still $f^S$ invariant. Here we recall a basic fact for recurrent times which is the claim in the proof of Proposition 3.1  \cite{KOTian}. This fact is that for any homeomorphism $g$, if $\nu$ is $g$-invariant and $\Gamma$ is a set with $\nu(\Gamma)>0,$ then for  $\nu$ a.e. $x\in \Gamma,$ $$\chi^*_{\Gamma}(x):=\lim_{n\rightarrow+\infty}\frac1n\sum_{j=0}^{n-1}\chi_{\Gamma}(g^{j}(x))>0$$(Remark that in general this fact holds only  for  $\nu$ a.e. $x\in \Gamma$ and then $x\in \cup_{n\in \mathbb{Z}}g^n(\Gamma)$ but maybe not hold for $\nu$ a.e. $x\in M$ since $\nu$ is just $g$-invariant). This guarantees that  for  $\mu$ a.e. $x\in A_\epsilon$, $\chi^*_{A_\epsilon}(x)>0 $ since $\mu$ is $f^S$ invariant and  $\mu(A_\epsilon)>0.$ Fix such a point $x\in A_\epsilon$. Let $$t_0(x)=0<t_1(x)<t_2(x)<\cdots
$$
to be the all positive times
such that $f^{t_i(x)S}(x)\in A_\epsilon.$ Then
$$\lim_{i\rightarrow+\infty}\frac{i}{t_i(x)}=\lim_{t_i(x)\rightarrow+\infty}\frac1{t_i(x)}
\sum_{j=0}^{t_i(x)-1}\chi_{A_\epsilon}(f^{jS}(x))=\chi^*_{A_\epsilon}(x)>0.$$
 Recall  
 the definition of $A_\epsilon$ and the second condition of GDS. Then $$\frac{\|Df^{S}|_{E(f^{t_j(x)S}(x))}\|}{ m(Df^{S}|_{F(f^{t_j(x)S}(x))})}\leq -\epsilon+\lambda^{-1}$$ and
 $$\frac{\|Df^{(t_{j+1}(x)-t_j(x)-1)S}|_{E(f^{t_j(x)S+S}(x))}\|}
 { m(Df^{(t_{j+1}(x)-t_j(x)-1)S}|_{F(f^{t_j(x)S+S}(x))})}\leq \lambda.$$ Hence, for $n=t_i(x),$ $$c_n(x)= \log\frac{\|Df^{nS}|_{E(x)}\|}{ m(Df^{nS}|_{F(x)})}\leq \sum_{j=0}^{i-1} \log\frac{\|Df^{S}|_{E(f^{t_j(x)S}(x))}\|}{ m(Df^{S}|_{F(f^{t_j(x)S}(x))})}$$
 $$\,\,\,+\sum_{j=0}^{i-1} \log\frac{\|Df^{(t_{j+1}(x)-t_j(x)-1)S}|_{E(f^{t_j(x)S+S}(x))}\|}
 { m(Df^{(t_{j+1}(x)-t_j(x)-1)S}|_{F(f^{t_j(x)S+S}(x))})}$$
$$ \leq i\log(-\epsilon+\lambda^{-1})+i\log\lambda=i\log (1-\epsilon\lambda).\,\,\,\,\,\,$$
Thus $$c(x)=\lim_{n\rightarrow +\infty} \frac1n c_n(x)\leq \lim_{i\rightarrow+\infty}\frac{i\log (1-\epsilon\lambda)}{t_i(x)}=\chi^*_{A_\epsilon}(x)\log (1-\epsilon\lambda)<0.$$
Remark that  $\chi^*_{A_\epsilon}(x)>0$ and the estimate inequality of $c_n(x)$ play the crucial roles.

Now we can use Lemma \ref{Lem:1} to get that there is $N > 0$
such that $a_N(x) < 0$ for every $x \in M.$ Recall that $a_N$ is a continuous function and $M$ is compact.
 So $t:=\max_{x\in M}\{a_N(x)\}$ exists and must be negative. If $\tau=e^t$, then $0<\tau<1$ and for any $x \in M,$
 $$\frac{\|Df^N|_{E(x)}\|}{ m(Df^N|_{F(x)})}=e^{a_N(x)}\leq e^t= \tau.$$ So the nontrivial global GDS $TM=E\oplus F$  is a nontrivial global $(N,\tau)$-dominated splitting.
We complete the proof.\qed

\medskip

\begin{Rem}\label{Rem:Main1} Note that the assumption of ergodicity of $\mu$ is not necessary due to the used claim in the proof of Proposition 3.1  \cite{KOTian}. Moreover, if we do not assume $f$ to be minimal in the proof, it is easy to see that for any  invariant (not necessarily ergodic) measure $\mu$, $\inf\frac1n\int_M a_n(x) d\mu<0$ if and only if  $\mu(\cup_{\epsilon>0}{A_\epsilon})>0.$
\end{Rem}

\begin{Rem}\label{Rem:Main2}
This proof implies a fact that  if $\inf\frac1n\int_M \log\frac{\|Df^n|_{E(x)}\|}{ m(Df^n|_{F(x)})} d\mu<0$  holds with respect to a continuous $Df$-invariant splitting $TM=E\oplus F$ for all ergodic invariant measure $\mu$, then $TM=E\oplus F$  is a global dominated splitting.
\end{Rem}




\section{Difference of GDS and Dominated splitting} \label{GDS-DS}
To further illustrate the new notion of GDS, we
construct a simple example which firstly appeared in \cite{ST}. This diffeomorphism has global GDS for $S=1,\lambda=1$ but does not have global dominated splitting.

\begin{Ex}\label{Ex:GDS-DS} Let $g$ be a $C^r (r\geq1)$ increasing function on $[0,1]$,
 satisfying: $$g(0)=0,\,\,g(1)=1,\,\,g'(0)=g'(1)=\frac{3-\sqrt[]{5}}2,\,\,g(\frac12)=\frac12,
 \,\,g'(\frac12)=\frac{3+\sqrt[]{5}}2\,\,\text{and}$$
 $$\frac{3-\sqrt[]{5}}2 \leq g'(x) \leq \frac{3+\sqrt[]{5}}2,\forall x\in[0,1],\,\,\,g(x)<x,\,\,\forall\,\,
 x\in(0,\frac12),\,\,g(x)>x,\,\,\forall\,\,x\in(\frac12,1) .$$
 \setlength{\unitlength}{1mm}
  \begin{figure}[!htbp]
  \begin{center}
  \begin{picture}(80,60)(0,0)
  \put(0,0){\scalebox{1}[1]{\includegraphics[0,0][30,40]{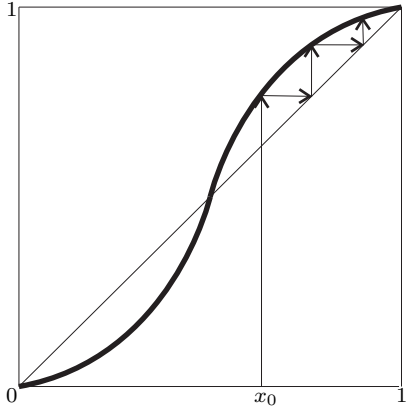}}}
  \put(6.6,54){$_{1}$}
\put(6.6,2.8){$_{0}$}\put(57.9,2.8){$_{1}$}\put(39.2,2.8){$_{x_0}$}
  \end{picture}
  \caption{Graph of the function $g$.}
  \label{pic:graph-of-g}
  \end{center}
  \end{figure}\\
And let $h: T^2\rightarrow T^2$ be the hyperbolic Torus automorphism
$$(y,z)\mapsto(2y+z,y+z),\,\, y,z\in \mathbb{S}^1=\mathbb{R}/\mathbb{Z}.
 $$ Define $f=g\times h:T^3\rightarrow
T^3$. Clearly,  $$Df(x,y,z)= \left(
\begin {array}{ccc}
g'(x)&0 &0\\
0&2&1\\
0&1&1
\end {array}
\right).
$$

There exists naturally a continuous splitting $TT^3=E_1\oplus
E_2\oplus E_3$, where $E_2$ and $E_3$ are from  the hyperbolic Torus
automorphism $h$ and $E_1$ is $g-$invariant. The forward  Lyapunov
exponent of $E_1$ is $\log\frac{3-\sqrt[]{5}}2$ over $T^3-\{\frac12\}\times T^2$
 and the Lyapunov
exponents of $ E_2\oplus E_3$ over $\{0\}\times T^2$ are
$\,\,\log\frac{3-\sqrt[]{5}}2,\,\,\log\frac{3+\sqrt[]{5}}2$
respectively. Denote by $\delta_0$ the point measure at point $0\in \mathbb{S}^1$ and
denote by $m$ the Lebesgue measure on $T^2$, then the product
measure $\mu=\delta_0\times m$ is a hyperbolic ergodic measure of
the diffeomorphism $f=g\times h$ with three nonzero Lyapunov exponents
$-\log\frac{3+\sqrt[]{5}}2,\,\,-\log\frac{3+\sqrt[]{5}}2,\,\,\log\frac{3+\sqrt[]{5}}2$. Set $E=E_1\oplus E_2$ and $F=E_3$, then $E\oplus
F$ construct a continuous $Df$-invariant splitting of $TT^3$ over
the whole space $T^3$ and $E\oplus
F$ is a
GDS  for $S=1,\lambda=1$ on the whole space $T^3$.
However, it is not a global dominated splitting since  for every point
$u=(\frac12,y,z)\,\,(y,z \in \mathbb{S}^1)$,
$$\frac{\|Df|_{E(u)}\|}{m(Df|_{F(u)})}=\frac{\frac{3+\sqrt[]{5}}2}{\frac{3+\sqrt[]{5}}2}=1.$$
Similarly  if $E=E_2$ and $F=E_1\oplus E_3$, we can follow above discussion to get that the new $E\oplus
F$ is also not dominated but a global GDS  for $S=1,\lambda=1$. But if $E=E_1$ and $F=E_2\oplus E_3$, then it is easy to see this splitting $E\oplus
F$ is not a  global GDS and thus is also not a global dominated splitting(even though this splitting is continuous on whole manifold).  All in all, this example has nontrivial global GDS but does not admit nontrivial global dominated splitting.
\qed
\end{Ex}

At the end of present paper, we point out a further question under a more general assumption.

\begin{Que}\label{Que:MainThm}
Let $f\in \Diff^1(M)$. If there is a nontrivial global  $Df-$invariant splitting $TM=E\oplus F$ satisfying \\
 (1). $TM=E\oplus F$  is continuous  on
$\Delta$;\\
(2). for every $x\in M,$ $$\liminf_{n\rightarrow +\infty}\frac1n\log\frac{\|Df^n|_{E(x)}\|}{m(Df^n|_{F(x)})}\leq0;$$ \\
 (3). there exists $x_0$, $$\liminf_{n\rightarrow +\infty}\frac1n\log\frac{\|Df^n|_{E(x_0)}\|}{m(Df^n|_{F(x_0)})}<0.$$ Then whether $f$ can not be minimal?
\end{Que}

{\bf Acknowledgement.} The author thanks very much
to Professor Marcus Bronzi, Krerley Oliveira  and Fernando Pereira Micena {\it et al} for their help in the period of  UFAL.

\bigskip
\end{document}